\documentclass{article}
\usepackage{amsmath}

\let\al\alpha
\let\be\beta
\let\ga\gamma
\let\de\delta
\let\eps\varepsilon
\let\ka\kappa
\let\la\lambda
\let\phi\varphi
\let\si\sigma

\newcommand\hu{\hat{u}}
\newcommand\oA{\bar{A}}
\newcommand\oE{\bar{E}}
\newcommand\oL{\bar{L}}
\newcommand\os{\bar{s}}
\newcommand\ou{\bar{u}}
\newcommand\ov{\bar{v}}
\newcommand\oy{\bar{y}}
\newcommand\oOm{\bar{\Omega}}
\newcommand\tchi{\tilde{\chi}}

\newcommand\cO{\mathcal{O}}
\newcommand\cX{\mathcal{X}}
\renewcommand\d{\mathrm{d}}

\renewcommand\Re{\operatorname{Re}}
\renewcommand\Im{\operatorname{Im}}

\begin{document}

\title{On the global error committed when evaluating the Evans
  function numerically%
  \footnote{This work was supported by EPSRC First Grant GR/S22134/01.}}
\author{Jitse Niesen\footnote{Mathematics Department, Heriot-Watt
    University, Edinburgh, EH14 4AS, United Kingdom. (Current address:
    Department of Mathematics, La Trobe University, Melbourne,
    Victoria 3086, Australia.) E-mail: j.niesen@latrobe.edu.au}}
\maketitle

\section{Introduction}

The Evans function is a tool for assessing the stability of travelling
waves solutions for partial differential equations.

A recent paper~\cite{malham.niesen:evaluating} analyzes the order
reduction experienced when evaluating the Evans function
numerically. The details of some lengthy calculations were excluded
from that paper for clarity. The purpose of this technical report is
to make these details publicly available. This report is not intended
to be read on its own; the reader is referred
to~\cite{malham.niesen:evaluating} for background and references.

\section{Setting}

We consider scalar reaction--diffusion equations of the form
$$
u_t = u_{xx} + f(u).
$$
Let $u(x,t) = \hu(\xi)$ with $\xi = x-ct$ be a travelling wave
solution of this equation. A linear stability analysis of this
travelling wave leads to the eigenvalue problem
\begin{subequations}
\label{bvp3}
\begin{equation}
\label{bvp3a} 
\frac{\d y}{\d\xi} = A(\xi;\la)\,y, 
\end{equation}
where
\begin{equation}
\label{bvp3b}
A(\xi;\la) = \begin{bmatrix}
0 & 1 \\ \la - f'(\hu(\xi)) & -c
\end{bmatrix}.
\end{equation}
\end{subequations}
The limits of~$A$ as $\xi\to\pm\infty$ are given by 
$$
A_\pm(\la) = \begin{bmatrix}
  0 & 1 \\ \la - f'(\hu_\pm) & -c 
\end{bmatrix}.
$$
Furthermore, the eigenvalues of~$A_-(\la)$ are 
\begin{align*}
\mu_-^{[1]} &= \frac12 \left( -c + \sqrt{c^2 + 4(\la-f'(\hu_-))}
\right) \quad\text{and} \\
\mu_-^{[2]} &= \frac12 \left( -c - \sqrt{c^2 + 4(\la-f'(\hu_-))}
\right),
\end{align*}
where $\hu_- = \lim_{\xi\to-\infty} \hu(\xi)$. Similarly, the
eigenvalues of~$A_+(\la)$ are  
\begin{align*}
\mu_+^{[1]} &= \frac12 \left( -c + \sqrt{c^2 + 4(\la-f'(\hu_+))}
\right) \quad\text{and} \\
\mu_+^{[2]} &= \frac12 \left( -c - \sqrt{c^2 + 4(\la-f'(\hu_+))}
\right),
\end{align*}
where $\hu_+ = \lim_{\xi\to+\infty} \hu(\xi)$.
The corresponding eigenvectors are $(1,\mu_\pm^{[1]})^\top$
and $(1,\mu_\pm^{[2]})^\top$, respectively.

To define the Evans function, assume that 
\begin{equation}
\label{lacond}
\Re \la > \max \big( f'(\hu_-), f'(\hu_+) \big) 
- \left( \frac{c}{\Im\la} \right)^2.
\end{equation}
This condition ensures that $\mu_-^{[1]}$ and $\mu_+^{[1]}$ have
positive real part, while $\mu_-^{[2]}$ and $\mu_-^{[2]}$ have
negative real part. 

The differential equation~\eqref{bvp3} is linear, and hence its
solutions form a linear space. Let $y_-$ be the solution which
satisfies 
\begin{equation}
\label{ymbc}
y_-(\xi) \sim \exp \big( \mu_-^{[1]}\xi \big)
\begin{bmatrix} 1 \\ \mu_-^{[1]} \end{bmatrix}
\qquad \text{as } \xi\to -\infty.
\end{equation}
Condition~\eqref{lacond} implies that this defines $y_-$ uniquely,
that $y_-$ satisfies the boundary condition $y(\xi) \to 0$
as $\xi \to -\infty$, and that any solutions satisfying this boundary
condition is a multiple of~$y_-$.

Similarly, we define $y_+$ as the solution satisfying
\begin{equation}
\label{ypbc}
y_+(\xi) \sim \exp \big( \mu_+^{[2]}\xi \big)
\begin{bmatrix} 1 \\ \mu_+^{[2]} \end{bmatrix}
\qquad \text{as } \xi\to \infty.
\end{equation}
The Evans function is then the function $D$ defined by
$$
D(\la) = \det \big[ \, y_-(0) \mid y_+(0) \, \big].
$$
We are interested in computing this function.

\section{Asymptotics near infinity}

In this section, we study the behaviour of $D(\la)$ as $|\la| \to \infty$. 

\subsection{The solution satisfying the left boundary condition}

Define a transformation $y_- \mapsto \oy_-$ by
\begin{equation}
\label{oym} 
y_-(\xi) = \exp(\mu_-^{[1]}\xi) \left( 
  \ou_-(\xi) \begin{bmatrix} 1 \\ \mu_-^{[1]} \end{bmatrix}
  + \ov_-(\xi) \begin{bmatrix} 1 \\ \mu_-^{[2]} \end{bmatrix}
\right)  = \exp(\mu_-^{[1]}\xi) \, B_- \, \oy_-(\xi),
\end{equation}
where
\begin{equation}
\oy_- = \begin{bmatrix} \ou_- \\ \ov_- \end{bmatrix} \text{ and } 
B_- = \begin{bmatrix} 1 & 1 \\ \mu_-^{[1]} & \mu_-^{[2]} \end{bmatrix}. 
\end{equation}
The differential equation~\eqref{bvp3} transforms to
\begin{equation}
\label{bvpoc}
\begin{aligned}
\frac{\d}{\d\xi} \ou_-(\xi) &= -\frac1{\ka_-} \phi_-(\xi) \ou_-(\xi) 
- \frac1{\ka_-} \phi_-(\xi) \ov_-(\xi), \\
\frac{\d}{\d\xi} \ov_-(\xi) &= \frac1{\ka_-} \phi_-(\xi) \ou_-(\xi) 
- \left( \ka_- - \frac1{\ka_-} \phi_-(\xi) \right) \ov_-(\xi),
\end{aligned}
\end{equation}
where
\begin{equation}
\label{ka} 
\phi_-(\xi) = f'(\hu(\xi)) - f'(\hu_-) \text{ and }
\ka_- = \sqrt{c^2 + 4(\la-f'(\hu_-))},
\end{equation}
and the boundary condition~\eqref{ymbc} becomes
\begin{equation}
\label{oymbc}
\ou_-(\xi) \to 1 \quad\text{and}\quad \ov_-(\xi) \to 0 \qquad\text{as }
\xi\to-\infty.
\end{equation}
Now, suppose that~$\ou_-$ and~$\ov_-$ can be expanded in
inverse powers of~$\ka_-$:
\begin{align*}
\ou_-(\xi;\ka_-) &= \ou_0^-(\xi) + \ka_-^{-1} \ou_1^-(\xi) + \ka_-^{-2} \ou_2^-(\xi) +
\ka_-^{-3} \ou_3^-(\xi) + \cO(\ka_-^{-4}), \\
\ov_-(\xi;\ka_-) &= \ov_0^-(\xi) + \ka_-^{-1} \ov_1^-(\xi) + \ka_-^{-2} \ov_2^-(\xi) +
\ka_-^{-3} \ov_3^-(\xi) + \cO(\ka_-^{-4}). 
\end{align*}
We now substitute these expansions in~\eqref{bvpoc} and equate the
coefficients of the powers of~$\ka_-$.
\begin{itemize}
\item At $\cO(\ka_-)$, we get $0 = \ov_0^-$, so $\ov_0^-$ is
identically zero.
\item At $\cO(1)$, we get $(\ou_0^-)' = 0$ and
$(\ov_0^-)' = -\ov_1^-$. The first equation, together with the boundary
condition~\eqref{oymbc}, implies that $\ou_0^- \equiv 1$. It follows
from the second equation that $\ov_1^-$  is identically zero. 
\item At $\cO(\ka_-^{-1})$, we get 
$$ 
(\ou_1^-)' = -\phi_-(\xi) \, (\ou_0^-+\ov_0^-) \quad\text{and}\quad
(\ov_1^-)' = \phi_-(\xi) \, (\ou_0^-+\ov_0^-) - \ov_2^-. 
$$
Substituting what we found before yields 
$$
(\ou_1^-)' = -\phi_-(\xi) \quad\text{and}\quad 0 = \phi_-(\xi) - \ov_2^-.
$$
Hence, $\ou_1^-(\xi) = -\Phi_-(\xi)$ and $\ov_2^-(\xi) = \phi_-(\xi)$,
where
$$
\Phi_-(\xi) = \int_{-\infty}^\xi \phi_-(x) \,\d{x}. 
$$
\item At $\cO(\ka_-^{-2})$, we get
$$ 
(\ou_2^-)' = -\phi_-(\xi) \, (\ou_1^-+\ov_1^-) \quad\text{and}\quad
(\ov_2^-)' = \phi_-(\xi) \, (\ou_1^-+\ov_1^-) - \ov_3^-. 
$$
Substituting $\ou_1^-$ and $\ov_1^-$ in the first equation yields
$(\ou_2^-)' = \phi_-(\xi) \, \Phi_-(\xi)$, and hence,
\begin{align*}
\ou_2^-(\xi) &= \int_{-\infty}^\xi \phi_-(x) \int_{-\infty}^x
\phi_-(y) \,\d{y} \,\d{x} \\
&= \frac12 \int_{-\infty}^\xi \int_{-\infty}^\xi \phi_-(x) \,
\phi_-(y) \,\d{y} \,\d{x} = \tfrac12 \big( \Phi_-(\xi) \big)^2.
\end{align*} 
\end{itemize}
Summarizing, we have the following series expansions for the
transformed solution:
\begin{equation}
\label{oy-exp}
\begin{aligned}
\ou_-(\xi;\ka_-) &= 1 - \ka_-^{-1} \Phi_-(\xi) + \tfrac12 \ka_-^{-2}
\big(\Phi_-(\xi)\big)^2 + \cO(\ka_-^{-3}), \\
\ov_-(\xi;\ka_-) &= \ka_-^{-2} \phi_-(\xi) + \cO(\ka_-^{-3}).
\end{aligned}
\end{equation}

\subsection{The solution satisfying the right boundary condition}

The computation for the solution~$y_+$ satisfying the right boundary
condition~\eqref{ypbc} is analogous. Define the transformation $y_+
\to \oy_+$ by 
$$
y_+(\xi) = \exp(\mu_+^{[2]}\xi) \left( 
  \ou_+(\xi) \begin{bmatrix} 1 \\ \mu_+^{[1]} \end{bmatrix}
  + \ov_+(\xi) \begin{bmatrix} 1 \\ \mu_+^{[2]} \end{bmatrix}
\right)  = \exp(\mu_+^{[2]}\xi) \, B_+ \, \oy_+(\xi) 
$$
where
$$
\oy_+ = \begin{bmatrix} \ou_+ \\ \ov_+ \end{bmatrix} \text{ and } 
B_+ = \begin{bmatrix} 1 & 1 \\ \mu_+^{[1]} & \mu_+^{[2]} \end{bmatrix}. 
$$
The differential equation~\eqref{bvp3} transforms to
\begin{equation}
\label{bvpocp}
\begin{aligned}
\frac{\d}{\d\xi} \ou_+(\xi) &= \left( \ka^+ - \frac1{\ka_+}
\phi_+(\xi) \right) \ou_+(\xi) - \frac1{\ka_+} \phi_+(\xi) \ov_+(\xi), \\
\frac{\d}{\d\xi} \ov_+(\xi) &= \frac1{\ka_+} \phi_+(\xi) \ou_+(\xi) 
+ \frac1{\ka_+} \phi_+(\xi) \ov_+(\xi),
\end{aligned}
\end{equation}
where
\begin{equation}
\label{kap} 
\phi_+(\xi) = f'(\hu(\xi)) - f'(\hu_+) \text{ and }
\ka_+ = \sqrt{c^2 + 4(\la-f'(\hu_+))},
\end{equation}
and the boundary condition~\eqref{ypbc} becomes
$$
\ou_+(\xi) \to 0 \quad\text{and}\quad \ov_+(\xi) \to 1 \qquad\text{as }
\xi\to\infty.
$$
Expand~$\ou_+$ and~$\ov_+$ in inverse powers of~$\ka_+$:
\begin{align*}
\ou_+(\xi;\ka_+) &= \ou_0^+(\xi) + \ka_+^{-1} \ou_1^+(\xi) + \ka_+^{-2} \ou_2^+(\xi) +
\ka_+^{-3} \ou_3^+(\xi) + \cO(\ka_+^{-4}), \\
\ov_+(\xi;\ka_-) &= \ov_0^+(\xi) + \ka_+^{-1} \ov_1^+(\xi) + \ka_+^{-2} \ov_2^+(\xi) +
\ka_+^{-3} \ov_3^+(\xi) + \cO(\ka_+^{-4}). 
\end{align*}
Substituting these expansions in the transformed differential
equation~\eqref{bvpocp} and equating the coefficients of the powers
of~$\ka_+$ yields the following equations:
\begin{align*}
0 &= \ou_0^+, & 0 &= 0, \\ 
(\ou_0^+)' &= -\ou_1^+, & (\ov_0^+)' &= 0,  \\
(\ou_1^+)' &= -\phi_+(\xi) \, (\ou_0^++\ov_0^+) + \ou_2^+, &
(\ov_1^+)' &= \phi_+(\xi) \, (\ou_0^++\ov_0^+), \\
(\ou_2^+)' &= -\phi_+(\xi) \, (\ou_1^++\ov_1^+) + \ou_3^+, &
(\ov_2^+)' &= \phi_+(\xi) \, (\ou_1^++\ov_1^+). 
\end{align*}
Solving these equations in the same way as in the previous case yields
the following series expansions for the transformed solution:
\begin{equation}
\label{oy-exp2}
\begin{aligned}
\ou_+(\xi;\ka_+) &= \ka_+^{-2} \phi_+(\xi) + \cO(\ka_+^{-3}), \\
\ov_+(\xi;\ka_+) &= 1 - \ka_+^{-1} \Phi_+(\xi) + \tfrac12 \ka_+^{-2}
\big(\Phi_+(\xi)\big)^2 + \cO(\ka_+^{-3}), 
\end{aligned}
\end{equation}
where
$$
\Phi_+(\xi) = \int_\xi^\infty \phi_+(x) \,\d{x}. 
$$

\subsection{Asymptotics for the Evans function}

The Evans function is obtained by evaluating both~$y_-$ and~$y_+$ at
$\xi=0$ and taking the wedge product (which equals the determinant of
the 2-by-2 matrix having these vectors as its columns). This yields
\begin{align*}
D(\la) &= y_=(0) \wedge y_+(0) 
= \big( B_-\oy_-(0) \big) \wedge \big( B_+\oy_+(0) \big) \\
&= 
\begin{bmatrix} 1 & 1 \\ \frac12(\ka_--c) & -\frac12(\ka_-+c) \end{bmatrix}
\begin{bmatrix} \ou_-(0) \\ \ov_-(0) \end{bmatrix} \wedge
\begin{bmatrix} 1 & 1 \\ \frac12(\ka_+-c) & -\frac12(\ka_++c) \end{bmatrix}
\begin{bmatrix} \ou_+(0) \\ \ov_+(0) \end{bmatrix} \\
&= \tfrac12 (\ka_--\ka_+) \big(\ov_-(0)\ov_+(0) - \ou_-(0)\ou_+(0) \big) \\
&\hspace{1cm} + \tfrac12 (\ka_-+\ka_+) \big(\ov_-(0)\ou_+(0) - \ou_-(0)\ov_+(0) \big)
\end{align*}
It follows from~\eqref{oy-exp} and~\eqref{oy-exp2} that, as $|\la|\to\infty$,
\begin{align*}
\ou_- &= \cO(1), & \ou_+ &= \cO(\ka_+^{-2}) = \cO(\la^{-1}), \\
\ov_- &= \cO(\ka_-^{-2}) = \cO(\la^{-1}), & \ov_+ &= \cO(1).
\end{align*}
Furthermore, we have $\ka_--\ka_+ = \cO\big(\la^{-1/2}\big)$. Hence,
$$
D(\la) = -\tfrac12 (\ka_- + \ka_+) \ou(0)\ov_+(0) + \cO(\la^{-1}).
$$
Finally, 
$$
D(\la) =  -2\la^{1/2} + \Phi - \tfrac14 \la^{-1/2} \Big(
\Phi^2 - 2f'(\hu_-) - 2f'(\hu_+) + c^2 \Big) + \cO(\la^{-1}), 
$$
where
$$
\Phi = \Phi_-(0) + \Phi_+(0) = \int_{-\infty}^0 \phi_-(x)
\,\d{x} + \int_0^\infty \phi_+(x) \,\d{x}. 
$$

\section{The exponential midpoint rule}

The exponential midpoint rule (or second-order Magnus method) for
solving the differential equation $y' = A(\xi) \, y$ is
$$
y_{k+1} = \exp\big( hA(\xi_k+\tfrac12h) \big) y_k.
$$
We assume that the step size~$h$ is fixed, so that $\xi_k = \xi_0 +
nh$. 

The transformation~\eqref{oym} changes the recursion for the
exponential midpoint rule to
\begin{align*}
\oy_{k+1} &= \exp(-\mu_-^{[1]}\xi_{n+1}) \, B_-^{-1} \, y_{k+1} \\
&= \exp(-\mu_-^{[1]}\xi_{k+1}) \, B_-^{-1} \, \exp\big( hA(\xi_k+\tfrac12h)
\big) y_k \\
&= \exp(-\mu_-^{[1]}\xi_{k+1}) \, B_-^{-1} \, \exp\big( hA(\xi_k+\tfrac12h)
\big) \exp(\mu_-^{[1]}\xi_k) \, B_- \, \oy_k \\
&= \exp(-h\mu_-^{[1]}) \, B_-^{-1} \, \exp\big( hB_-^{-1}A(\xi_k+\tfrac12h)B_-
\big) \exp(\mu_-^{[1]}\xi_k) \, \oy_k \\
&= \exp\big(h\oA_-(\xi_k+\tfrac12h) \big) \oy_k,
\end{align*}
where
$$
\oA_-(\xi) = B_-^{-1}A(\xi)B_- - \mu_-^{[1]} I.
$$
Using~\eqref{bvp3b}, we find that
$$
\oA_-(\xi) = \begin{bmatrix}
-\frac{1}{\ka_-} \phi_-(\xi) & -\frac{1}{\ka_-} \phi_-(\xi) \\[\jot]
\frac{1}{\ka_-} \phi_-(\xi) & -\ka_- + \frac{1}{\ka_-} \phi_-(\xi)
\end{bmatrix},
$$
with $\phi_-$ and $\ka_-$ as defined in~\eqref{ka}. 

The transformed recursion for the exponential midpoint rule is the
same as the exponential midpoint rule applied to the transformed
equation $\oy' = \oA(\xi)\,\oy$, cf.~\eqref{bvpoc}. The reason for this
is that the Magnus method is equivariant under transformations such
as~\eqref{oym}.

\subsection{The local error}

The local error for the exponential midpoint rule is defined by
$$
L^-_k = \exp\big(hA(\xi_k+\tfrac12h)\big) \, y(\xi_k) - y(\xi_{k+1}), 
$$
or, in transformed coordinates,
$$
\oL^-_k = \exp\big(h\oA_-(\xi_k+\tfrac12h)\big) \, \oy(\xi_k) - \oy(\xi_{k+1}), 
$$
We compute the matrix exponential of~$h\oA$ by diagonalization. The
eigenvalues of~$h\oA_-(\xi)$ are
\begin{align*}
\la_1 &= -\tfrac12 h \Bigl( \ka_- - \sqrt{\ka_-^2 - 4\phi_-(\xi)}
\Bigr) \\
&= -h \Bigl( \ka_-^{-1} \phi_-(\xi) + \ka^{-3} \big(\phi_-(\xi)\big)^2
+ \cO(\ka^{-5}) \Bigr)
\end{align*}
and
\begin{align*}
\la_2 &= -\tfrac12 h \Bigl( \ka_- + \sqrt{\ka_-^2 - 4\phi_-(\xi)}
\Bigr) \\
&= -h \Bigl( \ka_- - \ka_-^{-1} \phi_-(\xi) - \ka_-^{-3} \big(\phi_-(\xi)\big)^2
+ \cO(\ka_-^{-5}) \Bigr).
\end{align*}
The corresponding eigenvector matrix is
$$
V = \begin{bmatrix} 1 & 1 \\ v_1 & v_2 \end{bmatrix}
$$
with
\begin{align*}
v_1 &= \tfrac{1}{2\phi} \big( \ka^2 - \ka\sqrt{\ka^2-4\phi} \big) - 1, \\
v_2 &= \tfrac{1}{2\phi} \big( \ka^2 + \ka\sqrt{\ka^2-4\phi} \big) - 1,
\end{align*}
where we are writing $\ka$ for $\ka_-$ and $\phi$ for $\phi_-(\xi)$. 
Its inverse is
$$
V^{-1} = \frac{1}{v_2-v_1} 
\begin{bmatrix} v_2 & -1 \\ -v_1 & 1 \end{bmatrix}.
$$
We have $h\oA = V \big[ \begin{smallmatrix} \la_1 & 0 \\ 0 & \la_2
\end{smallmatrix} \big] V^{-1}$, and thus $\exp(h\oA) = V \big[
\begin{smallmatrix} \exp\la_1 & 0 \\ 0 & \exp\la_2 \end{smallmatrix}
\big] V^{-1}$. However, $\la_2 \sim -h\ka$ so that $\exp\la_2$ is
exponentially small as $|\ka|\to\infty$ if $\ka$ is restricted to lie
in a sector of the form $|\arg\ka| < \frac12\pi - \eps$. Under this
assumption,
\begin{align}
\exp(h\oA) &= V \begin{bmatrix} \exp\la_1 & 0 \\ 0 & \text{e.s.t.} 
\end{bmatrix} V^{-1} \notag\\
&= \frac{\exp\la_1}{v_2-v_1} \begin{bmatrix} v_2 & -1 \\ v_1v_2 & -v_1
\end{bmatrix} + \text{e.s.t.} \notag\\
&= \frac{\exp\bigl(-\frac12h\bigl(\ka-\sqrt{\ka^2-4\phi}\bigr)\bigr)}%
{\ka\sqrt{\ka^2-4\phi}} \notag\\
& \qquad {} \times\begin{bmatrix} 
\frac12\big(\ka^2+\ka\sqrt{\ka^2-4\phi}\big) - \phi & -\phi \\
\phi & -\frac12\big(\ka^2-\ka\sqrt{\ka^2-4\phi}\big) + \phi
\end{bmatrix} + \text{e.s.t.} \notag\\
\label{eOm-m2}
&= \begin{bmatrix}
1 - \dfrac{h\phi}{\ka} + \dfrac{h^2\phi^2}{2\ka^2} + \cO(\ka^{-3}) 
& -\dfrac{\phi}{\ka^2} + \cO(\ka^{-3}) \\[3\jot]
\dfrac{\phi}{\ka^2} + \cO(\ka^{-3}) &
-\dfrac{\phi^2}{\ka^4} + \cO(\ka^{-5})
\end{bmatrix},
\end{align}
where e.s.t. stands for exponentially small terms.

Hence, using the exact solution~\eqref{oy-exp}, we find that
\begin{multline*}
\exp\big(h\oA(\xi_k+\tfrac12h) \big) \oy(\xi_k) \\
{} = \begin{bmatrix}
1 - \dfrac{\Phi(\xi_k) + h\phi(\xi_k+\frac12h)}{\ka} 
+ \dfrac{\big(\Phi(\xi_k) + h\phi(\xi_k+\frac12h)\big)^2}{2\ka^2} +
\cO(\ka^{-3}) \\[3\jot]
\dfrac{\phi(\xi_k+\frac12h)}{\ka^2} 
+ \cO(\ka^{-3})
\end{bmatrix}
\end{multline*}
We thus arrive at the following expression for the local error:
\begin{align*}
\oL^-_k &= 
\exp\big(h\oA(\xi_k+\tfrac12h)\big) \, \oy(\xi_k) - \oy(\xi_{k+1}) \\
&= \begin{bmatrix}
\frac{\Phi(\xi_k+h) - (\Phi(\xi_k) + h\phi(\xi_k+\frac12h))}{\ka}
+ \frac{(\Phi(\xi_k+h))^2 
         - (\Phi(\xi_k) + h\phi(\xi_k+\frac12h))^2}{2\ka^2} +
\cO(\ka^{-3}) \\[\jot]
\frac{\phi(\xi_k+\frac12h) - \phi(\xi_{k+1})}{\ka^2} + \cO(\ka^{-3})
\end{bmatrix} \\
&= \begin{bmatrix}
\frac1{24} \ka^{-1} h^3 \phi''(\xi_k+\frac12h) +
\cO(\ka^{-1}h^5,\ka^{-2}h^3) \\[\jot]
-\frac12 \ka^{-2} h \phi'(\xi_k+\frac12h) +
\cO(\ka^{-2}h^2,\ka^{-3}h) 
\end{bmatrix}.
\end{align*}
We need to assume that $|\ka| \gg h^{-1}$ for the last equality.

\subsection{The global error}

We write the local error as
\begin{equation}
\label{oL-m2}
\oL^-_k = \begin{bmatrix}
\ka^{-1} \ga_k + \cO(\ka^{-2}h^3) \\
\ka^{-2} \de_k + \cO(\ka^{-3}h)
\end{bmatrix},
\end{equation}
where
\begin{align*}
\ga_k &= \int_{\xi_k}^{\xi_k+h} \phi_-(x) \,\d{x} - h
\phi_-(\xi_k+\tfrac12h) = \cO(h^3), \\
\de_k &= \phi_-(\xi_k+\tfrac12h) - \phi_-(\xi_k+h) = \cO(h).
\end{align*}
The global error satisfies the recurrence relation 
$$
E^-_{k+1} = \exp\bigl(hA(\xi_k+\tfrac12h)\bigr) \, E^-_k + L^-_k, \qquad E^-_0 = 0,
$$
or, in transformed coordinates,
$$
\oE^-_{k+1} = \exp\bigl(h\oA_-(\xi_k+\tfrac12h)\bigr) \, \oE^-_k + \oL^-_k, \qquad \oE^-_0 = 0.
$$
The solution of this recursion is
\begin{equation}
\label{oEk-ih}
\oE^-_k = \begin{bmatrix}
  \ka^{-1} \sum_{j=0}^{k-1} \ga_j + \cO(\ka^{-2}h^2) \\[\jot] 
  \ka^{-2} \de_{k-1} + \cO(\ka^{-3}h) 
\end{bmatrix}.
\end{equation}
This can easily be proved by induction. The case $k=1$ is trivial.
Assuming that \eqref{oEk-ih} holds for a particular value of~$k$, we
have, using~\eqref{eOm-m2} and~\eqref{oL-m2},
\begin{align*}
\oE^-_{k+1} &= \exp\bigl(h\oA_-(\xi_k+\tfrac12h)\bigr) \, \oE^-_k + \oL^-_k \\
&= \begin{bmatrix} 
1+\cO(\ka^{-1}) & \cO(\ka^{-2}) \\ \cO(\ka^{-2}) & \cO(\ka^{-4})
\end{bmatrix}
\begin{bmatrix}
  \ka^{-1} \sum_{j=0}^{k-1} \ga_j + \cO(\ka^{-2}h^2) \\[\jot] 
  \ka^{-2} \de_{k-1} + \cO(\ka^{-3}h) 
\end{bmatrix} \\
& \qquad\qquad +
\begin{bmatrix}
\ka^{-1} \ga_k + \cO(\ka^{-2}h^3) \\[\jot]
\ka^{-2} \de_k + \cO(\ka^{-3}h)
\end{bmatrix} \\
&= \begin{bmatrix}
  \ka^{-1} \sum_{j=0}^{k} \ga_j + \cO(\ka^{-2}h^2) \\[\jot] 
  \ka^{-2} \de_k + \cO(\ka^{-3}h) 
\end{bmatrix},
\end{align*}
which concludes the induction.

Substituting $\ga_k$ and $\de_k$ back in~\eqref{oEk-ih}, we find that
\begin{align*}
\oE^-_k &= \begin{bmatrix}
  \ka^{-1} \left( \int_{\xi_0}^{\xi_k} \phi_-(x) \,\d{x} -
     h \sum_{j=0}^{k-1} \phi_-(\xi_j+\tfrac12h) \right) + \cO(\ka^{-2}h^2) \\[2\jot]
  \ka^{-2} \big( \phi_-(\xi_k+\tfrac12h) - \phi_-(\xi_k+h) \big)
     + \cO(\ka^{-3}h) 
\end{bmatrix} \\
&= \begin{bmatrix}
  \ka^{-1} \left( \int_{\xi_0}^{\xi_k} \phi(x) \,\d{x} -
     h \sum_{j=0}^{k-1} \phi(\xi_j+\tfrac12h) \right) + \cO(\ka^{-2}h^2) \\[2\jot]
  \ka^{-2} \big( \phi(\xi_k+\tfrac12h) - \phi(\xi_k+h) \big)
     + \cO(\ka^{-3}h) 
\end{bmatrix}
\end{align*}
where $\phi(\xi) = f'(\hu(\xi))$ differs from $\phi_-(\xi)$ by a constant.

\subsection{The solution on $[0,\infty)$}

To compute the solution~$y_+$ satisfying the right boundary condition,
we run the exponential midpoint rule backwards:
$$
\oy_{k+1} = \exp\bigl(-h\oA_+(\xi_k-\tfrac12h) \bigr) \, \oy_k,
$$
where $\xi_k = \xi_0 - kh$ and
\begin{equation}
\label{oAp}
\oA_+(\xi) = B_+^{-1}A(\xi)B_+ - \mu_+^{[2]} I
= \begin{bmatrix}
\ka_+ - \frac{1}{\ka_+} \phi_+(\xi) & -\frac{1}{\ka_+} \phi_+(\xi) \\[\jot]
\frac{1}{\ka_+} \phi_+(\xi) & \frac{1}{\ka_+} \phi_+(\xi)
\end{bmatrix},
\end{equation}
with $\ka_+$ and $\phi_+$ as defined in~\eqref{kap}. A similar
computation as before yields 
$$
\exp(-h\oA_+) = \begin{bmatrix}
-\dfrac{\phi_+^2}{\ka_+^4} + \cO(\ka_+^{-5}) & 
\dfrac{\phi_+}{\ka_+^2} + \cO(\ka_+^{-3}) \\[3\jot]
-\dfrac{\phi_+}{\ka_+^2} + \cO(\ka_+^{-3}) &
1 - \dfrac{h\phi_+}{\ka_+} + \dfrac{h^2\phi_+^2}{2\ka_+^2} + \cO(\ka_+^{-3}) 
\end{bmatrix}
$$
and
\begin{multline*}
\exp\big(h\oA_+(\xi_k-\tfrac12h) \big) \, \oy(\xi_k) \\
{} = \begin{bmatrix}
\dfrac{\phi_+(\ldots)}{\ka_+^2} + \cO(\ka_+^{-3}) \\[3\jot]
1 - \dfrac{\Phi_+(\xi_k) + h\phi_+(\ldots)}{\ka_+} 
+ \dfrac{\big(\Phi_+(\xi_k) + h\phi_+(\ldots)\big)^2}{2\ka_+^2} +
\cO(\ka_+^{-3}) 
\end{bmatrix}
\end{multline*}
where $\phi_+(\ldots)$ stands for $\phi_+(\xi_k-\frac12h)$. When we
use this to determine the local error, we find
$$
\oL^+_k
= \begin{bmatrix}
\dfrac{\phi_+(\xi_k-\tfrac12h) - \phi_+(\xi_{k+1})}{\ka_+^2} 
+ \cO(\ka_+^{-3}) \\[3\jot]
\dfrac{\Phi_+(\xi_{k+1}) - \Phi_+(\xi_k) - h\phi_+(\xi_k+\frac12h)}{\ka_+}
+ \cO(\ka_+^{-2})
\end{bmatrix}.
$$
A standard induction argument shows that the global error is
$$
\oE^+_k
= \begin{bmatrix}
  \ka_+^{-2} \big( \phi(\xi_k+\tfrac12h) - \phi(\xi_k+h) \big) 
     + \cO(\ka_+^{-3}h) \\[2\jot]
  \ka_+^{-1} \left( \int_{\xi_k}^{\xi_0} \phi(x) \,\d{x} -
     h \sum_{j=0}^{k-1} \phi(\xi_j-\tfrac12h) \right) + \cO(\ka_+^{-2}h^2) 
\end{bmatrix}.
$$

\subsection{The error in the Evans function}

The numerically computed value for the Evans function is the wedge
product of the numerical solutions:
$$
D_{\text{num}}(\la) = y_k^- \wedge y_k^+ = 
\bigl( y_-(0) + E^-_k \bigr) \wedge \bigl( y_+(0) + E^+_k \bigr),
$$
with $k$ chosen such that $0=0$. Hence the error in the Evans
function is
\begin{align*}
E_D(\la) &= D(\la) - D_{\text{num}}(\la) \\
&= y_-(0) \wedge E^+_k + E^-_k \wedge y_+(0) 
+ E^-_k \wedge E^+_k \\
&= B_-\oy_-(0) \wedge B_+\oE^+_k + B_-\oE^-_k \wedge B_+\oy_+(0) 
+ B_-\oE^-_k \wedge B_+\oE^+_k. 
\end{align*}
Substituting $B_-$ and $B_+$, we find
\begin{equation}
\label{ed}
\begin{aligned}
E_D &= \tfrac12(\ka_- - \ka_+) \Bigl( 
\ov_-(0) \, [\oE^+_k]_2 - \ou_-(0) \, [\oE^+_k]_1 
+ [\oE^-_k]_2 \, \ov_+(0) \\ 
&\hspace{20ex} - [\oE^-_k]_1 \, \ou_+(0) 
+ [\oE^-_k]_2 \, [\oE^+_k]_2 - [\oE^-_k]_1 \, [\oE^+_k]_1
\Bigr) \\
&\quad + \tfrac12(\ka_- + \ka_+) \Bigl( 
\ov_-(0) \, [\oE^+_k]_1 - \ou_-(0) \, [\oE^+_k]_2 
+ [\oE^-_k]_2 \, \ou_+(0) \\ 
&\hspace{20ex} - [\oE^-_k]_1 \, \ov_+(0) 
+ [\oE^-_k]_2 \, [\oE^+_k]_1 - [\oE^-_k]_1 \, [\oE^+_k]_2
\Bigr).
\end{aligned}
\end{equation}
We now estimate all the terms in this expression and drop the ones of
lower order:
\begin{align*}
E_D &= - \tfrac12(\ka_- + \ka_+) \Bigl( \ou_-(0) \, [\oE^+_k]_2 
+ [\oE^-_k]_1 \, \ov_+(0) + [\oE^-_k]_1 \, [\oE^+_k]_2 \Bigr) +
\cO(\la^{-1}h) \\
&= h \sum_{j=-N}^{N-1} \phi(jh+\tfrac12h) - \int_{-L}^L \phi(x) \,\d{x}
+ \cO(\la^{-1/2}h^2),
\end{align*}
assuming that the differential equations are solved on the
intervals~$[-L,0]$ and~$[0,L]$ with $L=Nh$. The final step is to apply
the Euler--MacLaurin summation formula (see
e.g.~\cite{atkinson:introduction}), which states that
\begin{align*}
h \sum_{j=0}^n f(jh) &= \int_0^{nh} f(x) \,\d{x} 
+ \tfrac12 h \bigl( f(0)+f(nh) \bigr) \\
&\qquad\qquad + \sum_{i=1}^m \frac{B_{2i}h^{2i}}{(2i)!} 
\bigl( f^{(2i-1)}(nh) - f^{(2i-1)}(0) \bigr) \\
&\qquad\qquad + \frac{nB_{2m+2}h^{2m+3}}{(2m+2)!} f^{(2m+2)}(\xi)
\end{align*}
for some $\xi \in [0,nh]$, where $B_k$ denote the Bernoulli
numbers. This yields
\begin{multline*}
E_D = \int_{-L}^{-L+\frac12h} \phi(x) \,\d{x} 
+ \int_{L-\frac12h}^L \phi(x) \,\d{x} 
+ \tfrac12 h \bigl( \phi(-L) + \phi(L) \bigr) \\
\qquad\qquad + \sum_{i=1}^m \frac{B_{2i}h^{2i}}{(2i)!} 
\bigl( f^{(2i-1)}(L) - f^{(2i-1)}(-L) \bigr) 
+ \cO(h^{2m+2}, \la^{-1/2}h^2).
\end{multline*}
Now, $\phi(\xi)$ decays exponentially fast to zero as $|\xi| \to
\infty$. So if we assume that $L$ is sufficiently large, we can ignore
all terms in this equation but the last one. We thus arrive at the
final result, which is that the error in the Evans function is of
order $\la^{-1/2}h^2$.

\section{The fourth-order Magnus method}

We repeat the computation in the previous section for the fourth-order
Magnus method. This method is given by
$$
y_{k+1} = \exp\Bigl( 
\tfrac12h \bigl( A(\xi_k^1) + A(\xi_k^2) \bigr)
- \tfrac{\sqrt3}{12} h^2 \bigl[ A(\xi_k^1), A(\xi_k^2) \bigr] 
\Bigr) y_k,
$$
where $[\,\cdot\,,\,\cdot\,]$ denotes the matrix commutator defined by
$[X,Y] = XY - YX$ and $\xi_k^1$, $\xi_k^2$ are the Gauss--Legendre
points 
\begin{equation}
\label{gl-points}
\xi_k^1 = \xi_k + (\tfrac12-\tfrac16\sqrt3) h \quad\text{and}\quad
\xi_k^2 = \xi_k + (\tfrac12+\tfrac16\sqrt3) h.
\end{equation}
After the transformation~\eqref{oym}, the method reads
$$
\oy_{k+1} = \exp(\oOm_k) \oy_k
$$
with
\begin{align*}
\oOm_k &= \tfrac12h \bigl( \oA_-(\xi_k^1) + \oA_-(\xi_k^2) \bigr)
- \tfrac{\sqrt3}{12} h^2 \bigl[ \oA_-(\xi_k^1), \oA_-(\xi_k^2) \bigr] 
\\
&= h \begin{bmatrix} 
  -\dfrac{\al_k}{\ka_-} & \be_k - \dfrac{\al_k}{\ka_-} \\[3\jot]
  \be_k + \dfrac{\al_k}{\ka_-} & -\ka_- + \dfrac{\al_k}{\ka_-}
\end{bmatrix},
\end{align*}
where $\al_k$ and $\be_k$ are given by
\begin{equation}
  \label{albe}
  \al_k = \tfrac12 \bigl( \phi_-(\xi_k^1)+\phi_-(\xi_k^2) \bigr)
  \quad\text{and}\quad
  \be_k = -\frac{\sqrt{3}}{12} h \bigl(
  \phi_-(\xi_k^1)-\phi_-(\xi_k^2) \bigr). 
\end{equation}

\subsection{The local error}

The eigenvalues of~$\oOm_k$ are
$$
\la_1 = -\tfrac12 h (\ka_- - \tchi_k) 
= -h \Bigl( \ka_-^{-1} \chi_k + \ka^{-3} \chi_k^2 + \cO(\ka^{-5}) \Bigr)
$$
and
$$
\la_2 = -\tfrac12 h (\ka_- + \tchi_k)
= -h \Bigl( \ka_- - \ka_-^{-1} \chi_k - \ka_-^{-3} \chi_k^2 + \cO(\ka_-^{-5}) \Bigr),
$$
where $\chi_k = \al_k - \be_k^2$ and
$$
\tchi_k = \sqrt{\ka_-^2 - 4(\al_k - \be_k^2)}.
$$
The corresponding eigenvector matrix is
$$
V = \begin{bmatrix} 1 & 1 \\ v_1 & v_2 \end{bmatrix}
$$
with
$$
v_1 = 
\frac{2\al_k - \ka_-^2 + \ka_-\tchi_k}{2(\ka_-\be_k + \al_k)}
\quad\text{and}\quad
v_2 = 
\frac{2\al_k - \ka_-^2 - \ka_-\tchi_k}{2(\ka_-\be_k + \al_k)}.
$$
Its inverse is
$$
V^{-1} = \frac{1}{v_2-v_1} 
\begin{bmatrix} v_2 & -1 \\ -v_1 & 1 \end{bmatrix}.
$$
As in the previous section, where we were considering the exponential
midpoint rule, we have $\la_1 \sim -h\ka$ so that $\exp\la_1$ is
exponentially small (under the same assumption as before). Hence,
\begin{align}
\exp(\oOm_k) &= V \begin{bmatrix} \exp\la_1 & 0 \\ 0 & \text{e.s.t.} 
\end{bmatrix} V^{-1} \notag\\
&= \frac{\exp\la_1}{v_2-v_1} \begin{bmatrix} v_2 & -1 \\ v_1v_2 & -v_1
\end{bmatrix} + \text{e.s.t.} \notag\\
&= \exp\bigl( -\tfrac12 h ( \ka_- - \tchi_k ) \bigr)
   \frac{\al_k - \ka_-\be_k}{\ka_-\tchi_k} \notag\\
& \qquad {} \times\begin{bmatrix} 
\dfrac{2\al_k - \ka_-^2 - \ka_-\tchi_k}{2(\ka_-\be_k-\al_k)} & -1 \\[3\jot]
-\dfrac{\ka_-\be_k + \al_k}{\ka_-\be_k - \al_k}
& -\dfrac{2\al_k - \ka_-^2 + \ka_-\tchi_k}{2(\ka_-\be_k-\al_k)}
\end{bmatrix} + \text{e.s.t.} \notag\\
&= \frac{\exp\bigl( -\tfrac12 h (\ka_- - \tchi_k) \bigr)}{\ka_-\tchi_k} 
\begin{bmatrix}
\frac12 (\ka_-^2 + \ka_-\tchi_k) - \al_k & \ka_-\be_k - \al_k \\[1\jot]
\ka_-\be_k + \al_k & -\frac12 (\ka_-^2 - \ka_-\tchi_k) - \al_k
\end{bmatrix} \notag \\
\label{eOm-m4}
&= \begin{bmatrix}
1 - \dfrac{h\chi_k}{\ka_-} + \dfrac{h^2\chi_k^2-2\be_k^2}{2\ka_-^2} 
+ \cO(\ka_-^{-3}) 
& \dfrac{\be_k}{\ka_-} + \cO(\ka_-^{-2}) \\[3\jot]
\dfrac{\be_k}{\ka_-} + \cO(\ka_-^{-2}) 
& \dfrac{\be_k}{\ka_-^2} + \cO(\ka_-^{-3})
\end{bmatrix},
\end{align}
where e.s.t. stands for exponentially small terms.

Hence, using the exact solution~\eqref{oy-exp}, we find that
$$
\exp(\oOm_k) \, \oy(\xi_k) = \begin{bmatrix}
1 - \dfrac{\Phi_-(\xi_k) + h\chi_k}{\ka_-} + \cO(\ka_-^{-2}) \\[3\jot]
\dfrac{\be_k}{\ka_-} + \cO(\ka_-^{-2})
\end{bmatrix}.
$$
We thus arrive at the following expression for the local error: 
$$
\oL^-_k =  \exp(\oOm_k) \, \oy(\xi_k) - \oy(\xi_{k+1}) 
= \begin{bmatrix} 
\ka_-^{-1} \ga_k + \cO(\ka_-^{-2}h^4) \\[\jot]
\ka_-^{-1} \be_k + \cO(\ka_-^{-2}h) 
\end{bmatrix}, 
$$
where
\begin{align*}
\ga_k &= 
\int_{\xi_k}^{\xi_{k+1}} \phi_-(x) \,\d{x} - h(\al_k-\be_k^2) \\
&= h^5 \left( \tfrac{1}{4320} \phi''''(\xi_k+\tfrac12h) + 
\tfrac1{144} \big(\phi'(\xi_k+\tfrac12h)\big)^2 \right) + \cO(h^7) 
\end{align*}
and
$$
\be_k = \tfrac1{12} h^2 \phi'(\xi_k+\tfrac12h) + \cO(h^4).
$$

\subsection{The global error}

The global error satisfies the recurrence relation 
$$
\oE^-_{k+1} = \exp(\oOm_k) \, \oE^-_k + \oL^-_k, \qquad \oE^-_0 = 0.
$$
The solution of this recursion is
$$
\oE^-_k = \begin{bmatrix}
  \ka_-^{-1} \sum_{j=0}^{k-1} \ga_j + \cO(\ka_-^{-2}h^4) \\[\jot] 
  \ka_-^{-1} \be_{k-1} + \cO(\ka_-^{-2}h) 
\end{bmatrix}.
$$
This can easily be proved by induction. The key step in the proof is
the following computation:
\begin{align*}
\oE^-_{k+1} &= \exp(\oOm_kA) \, \oE^-_k + \oL^-_k \\
&= \begin{bmatrix} 
1+\cO(\ka_-^{-1}) & \cO(\ka_-^{-1}) \\[\jot]
 \cO(\ka_-^{-1}) & \cO(\ka_-^{-2})
\end{bmatrix}
\begin{bmatrix}
  \ka_-^{-1} \sum_{j=0}^{k-1} \ga_j + \cO(\ka_-^{-2}h^4) \\[\jot] 
  \ka_-^{-1} \be_{k-1} + \cO(\ka_-^{-2}h) 
\end{bmatrix} \\
& \qquad\qquad +
\begin{bmatrix}
\ka_-^{-1} \ga_k + \cO(\ka_-^{-2}h^4) \\[\jot]
\ka_-^{-1} \be_k + \cO(\ka_-^{-2}h) 
\end{bmatrix} \\
&= \begin{bmatrix}
  \ka_-^{-1} \sum_{j=0}^k \ga_j + \cO(\ka_-^{-2}h^4) \\[\jot] 
  \ka_-^{-1} \be_k + \cO(\ka_-^{-2}h) 
\end{bmatrix}.
\end{align*}
Substituting $\ga_j$ back in the formula for~$\oE^-_k$, we find that
$$
\oE^-_k = \begin{bmatrix}
  \ka^{-1} \left( \int_{\xi_0}^{\xi_k} \phi_-(x) \,\d{x} -
     h \sum_{j=0}^{k-1} (\al_k-\be_k^2) \right) + \cO(\ka^{-2}h^4) \\[2\jot]
  \ka^{-1} \be_k + \cO(\ka_-^{-2}h) 
\end{bmatrix}.
$$
We now approximate
$$
  \be_k = -\frac{\sqrt{3}}{12} h \bigl(
  \phi_-(\xi_k^1)-\phi_-(\xi_k^2) \bigr)
  = \tfrac1{12} h^2 \phi'(\xi_k-\tfrac12h) + \cO(h^4).
$$
Thus, 
$$
h \sum_{j=0}^{k-1} \be_k^2 = \tfrac1{12} h^4 \int_{\xi_0}^{\xi_k}
\bigl(\phi'(x)\bigr)^2 \,\d{x} + \cO(h^6),
$$
and therefore,
$$
\oE^-_k = \begin{bmatrix}
  \displaystyle \ka_-^{-1}\biggl( 
  \int_{\xi_0}^{\xi_k} \phi(x) + \tfrac1{144} h^4 (\phi'(x))^2 \,\d{x} 
	- \tfrac12 h \sum_{j=0}^{k-1} \bigl(\phi(\xi_k^1)+\phi(\xi_k^2)\bigr) \biggr)
  \hspace*{3ex} \\
  \hspace*{46ex} {} + \cO(\ka_-^{-1}h^6,\ka_-^{-2}h^4) \\[2\jot] 
  \frac1{12} \ka_-^{-1}h^2 \phi'(\xi_k-\frac12h) +
  \cO(\ka_-^{-1}h^4,\ka_-^{-2}h)
\end{bmatrix}.
$$

\subsection{The solution on $[0,\infty)$}

To compute the solution~$y_+$ satisfying the right boundary condition,
we run the same method backwards:
$$
\oy_{k+1} = \exp(\oOm^+_k) \, \oy_k
$$
with 
$$
\oOm^+_k = -\tfrac12h \bigl( \oA_+(\xi_k^1) + \oA_+(\xi_k^2) \bigr)
- \tfrac{\sqrt3}{12} h^2 \bigl[ \oA_+(\xi_k^1), \oA_+(\xi_k^2) \bigr],
$$
where
\begin{equation}
\label{xikp}
\xi_k = \xi_0 - kh, \qquad
\xi_k^1 = \xi_k - (\tfrac12-\tfrac16\sqrt3) h \quad\text{and}\quad
\xi_k^2 = \xi_k - (\tfrac12+\tfrac16\sqrt3) h.
\end{equation}
We can write the matrix $\oOm^+_k$ as
$$
\oOm^+_k = h \begin{bmatrix} 
  \ka_+ - \dfrac{\al_k}{\ka_+} & \be_k + \dfrac{\al_k}{\ka_+}
  \\[3\jot]
  \be_k - \dfrac{\al_k}{\ka_+} & - \dfrac{\al_k}{\ka_+}
\end{bmatrix},
$$
where
$$
  \al_k = \tfrac12 \bigl( \phi_+(\xi_k^1)+\phi_+(\xi_k^2) \bigr)
  \quad\text{and}\quad
  \be_k = -\frac{\sqrt{3}}{12} h \bigl(
  \phi_+(\xi_k^1)-\phi_+(\xi_k^2) \bigr). 
$$
A similar computation as before yields 
$$
\exp(-\oOm_+) = \begin{bmatrix}
  \dfrac{\be_k^2}{\ka_+^2} + \cO(\ka_+^{-3}) &
  \dfrac{\be_k}{\ka_+} + \cO(\ka_+^{-2}) \\[3\jot]
  \dfrac{\be_k}{\ka_+} + \cO(\ka_+^{-2}) &
  1 - \dfrac{h(\al_k-\be_k^2)}{\ka_+} +
  \dfrac{h^2(\al_k-\be_k^2)^2-2\be_k^2}{2\ka_+^2} + \cO(\ka_+^{-3})
\end{bmatrix}
$$
and
$$
\exp(-\oOm_+) \big) \, \oy(\xi_k) = \begin{bmatrix}
  \dfrac{\be_k}{\ka_+} + \cO(\ka_+^{-2}) \\[3\jot]
  1 - \dfrac{\Phi_+(\xi_k)+h(\al_k-\be_k^2)}{\ka_+} + \cO(\ka_+^{-2})
\end{bmatrix}.
$$
When we use this to determine the local error, we find
$$
\oL^+_k
= \begin{bmatrix}
  \dfrac{\be_k}{\ka_+} + \cO(\ka_+^{-2}h) \\
  \displaystyle \dfrac{1}{\ka_+} \biggl( \int_{\xi_{k+1}}^{\xi_k}
    \phi_+(x)\,\d{x} - h(\al_k-\be_k^2) \biggr) + \cO(\ka_+^{-2}h^4)
\end{bmatrix}.
$$
As in the previous section, we conclude that the global error is given
by
$$
\oE^+_k = \begin{bmatrix}
  \frac1{12} \ka_-^{-1}h^2 \phi'(\xi_k+\frac12h) +
  \cO(\ka_-^{-1}h^4,\ka_-^{-2}h) \\[4\jot]
  \displaystyle \ka_-^{-1}\biggl( 
  \int_{\xi_0}^{\xi_k} \phi(x) + \tfrac1{144} h^4 (\phi'(x))^2 \,\d{x} 
	- \tfrac12 h \sum_{j=0}^{k-1} \bigl(\phi(\xi_k^1)+\phi(\xi_k^2)\bigr) \biggr)
  \hspace*{3ex} \\
  \hspace*{46ex} {} + \cO(\ka_-^{-1}h^6,\ka_-^{-2}h^4)
\end{bmatrix}.
$$

\subsection{The error in the Evans function}

Substituting the results for the global error in~\eqref{ed}, we find
that the error in evaluating the Evans function is
\begin{multline*}
E_D = \tfrac12 h \sum_{j=-N}^{N-1} \Bigl(
\phi\bigl(jh+(\tfrac12-\tfrac16\sqrt3)h\bigr) +
\phi\bigl(jh+(\tfrac12+\tfrac16\sqrt3)h\bigr) \Bigr) \\
- \int_{-L}^L \phi(x) + \tfrac1{144} h^4 (\phi'(x))^2 \,\d{x} + \cO(\la^{-1/2}h^2),
\end{multline*}
As with the exponential midpoint rule, the Euler--MacLaurin summation
formula can be applied to show that the term
\begin{equation}
\label{olka}
\tfrac12 h \sum_{j=-N}^{N-1} \Bigl(
\phi\bigl(jh+(\tfrac12-\tfrac16\sqrt3)h\bigr) +
\phi\bigl(jh+(\tfrac12+\tfrac16\sqrt3)h\bigr) \Bigr) 
- \int_{-L}^L \phi(x) \,\d{x}
\end{equation}
is negligible if $L$ is sufficiently large. So, we find that the error
in the Evans function is given by
$$
E_D = - \tfrac1{144} h^4 \int_{-\infty}^\infty (\phi'(x))^2 \,\d{x} +
\cO(\la^{-1/2}h^2).
$$

\section{The fourth-order Gauss--Legendre method}

The two-stage Gauss--Legendre method for solving the equation $y' =
A(\xi)\,y$ is given by
\begin{align*}
s_1 &= A(\xi_k^1) \, 
\bigl( y_k + \tfrac14h s_1 + (\tfrac14-\tfrac{\sqrt3}6)h s_2 \bigr), \\
s_2 &= A(\xi_k^2) \, 
\bigl( y_k + (\tfrac14+\tfrac{\sqrt3}6)h s_1 + \tfrac14h s_2 \bigr), \\
y_{k+1} &= y_k + \tfrac12h (s_1 + s_2),
\end{align*}
where $\xi_k^1$ and $\xi_k^2$ are the Gauss--Legendre points, given
in~\eqref{gl-points}. As usual, we transform this to
\begin{align*}
\os_1 &= \oA_-(\xi_k^1) \, 
\bigl( \oy_k + \tfrac14h \os_1 + (\tfrac14-\tfrac{\sqrt3}6)h \os_2 \bigr), \\
\os_2 &= \oA_-(\xi_k^2) \, 
\bigl( \oy_k + (\tfrac14+\tfrac{\sqrt3}6)h \os_1 + \tfrac14h \os_2 \bigr), \\
\oy_{k+1} &= \oy_k + \tfrac12h (\os_1 + \os_2).
\end{align*}

\subsection{The local error}

Substituting $\oA_-$ and $\oy_k = \oy(\xi_k)$ in the above formula and
rearranging yields
\begin{align*}
& \left( 1+\frac{h\phi_1}{4\ka} \right) \os_{11} +
\frac{h\phi_1}{4\ka} \os_{12} + \frac{h\si_1\phi_1}{\ka} \os_{21} +
\frac{h\si_1\phi_1}{\ka} \os_{22} 
= -\frac{\phi_1}{\ka} \oy_{k1} - \frac{\phi_1}{\ka} \oy_{k2},
\\[4\jot]
& -\frac{h\phi_1}{4\ka} \os_{11} + \left( \frac{h\ka}4 + 1
-\frac{h\phi_1}{4\ka} \right) \os_{12} - \frac{h\si_1\phi_1}{\ka}
\os_{21} + h\si_1 \left( \ka - \frac{\phi_1}{\ka} \right) \os_{22}
\\
& \hspace{48ex} {} = \frac{\phi_1}{\ka} \oy_{k1} 
- \left( \ka - \frac{\phi_1}{\ka} \right) \oy_{k2},
\\[4\jot]
& \frac{h\si_2\phi_2}{\ka} \os_{11} + \frac{h\si_2\phi_2}{\ka}
\os_{12} + \left( 1+\frac{h\phi_2}{4\ka} \right) \os_{21} +
\frac{h\phi_2}{4\ka} \os_{22} = \oy_{k1}
= -\frac{\phi_2}{\ka} \oy_{k1} - \frac{\phi_2}{\ka} \oy_{k2},
\\[4\jot]
& - \frac{h\si_2\phi_2}{\ka} \os_{11} + h\si_2 \left( \ka -
\frac{\phi_2}{\ka} \right) \os_{12} - \frac{h\phi_2}{4\ka} \os_{21} +
\left( \frac{h\ka}4 + 1 - \frac{h\phi_2}{4\ka} \right) \os_{22}
\\
& \hspace{48ex} {} = \frac{\phi_2}{\ka} \oy_{k1} 
- \left( \ka - \frac{\phi_2}{\ka} \right) \oy_{k2},
\end{align*}
where 
$$
\si_1 = \tfrac14 - \tfrac16\sqrt3, \quad
\si_2 = \tfrac14 + \tfrac16\sqrt3, \quad 
\phi_1 = \phi_-(\xi_k^1) \quad\text{and}\quad
\phi_2 = \phi_-(\xi_k^2).
$$
We assume that the unknowns $\os_{ij}$ can be expanded in powers of $\ka$ like
$$
\os_{ij} = \os_{ij}^0 + \frac{\os_{ij}^1}{\ka} +
\frac{\os_{ij}^2}{\ka^2} + \cO(\ka^{-3}).
$$
Substitute this in the set of four equations above and collect like
powers of $\ka$. 
\begin{itemize}
\item At order~$\ka$, we find 
\begin{gather*}
h\os_{12}^0 + h\si_1\os_{22}^0 = -\oy_{k2}, \\
h\si_2\os_{12}^0 + h\os_{22}^0 = -\oy_{k2}.
\end{gather*}
Solving these equations yields
$$
\os_{12}^0 = -\frac{2\sqrt3}{h} \, \oy_{k2} \quad\text{and}\quad
\os_{22}^0 = \frac{2\sqrt3}{h} \, \oy_{k2}.
$$
\item At order~$\ka^0$, we find
\begin{gather*}
\os_{11}^0 = 0, \\
\tfrac14 h \os_{12}^1 + \os_{12}^0 + h \si_1 \os_{22}^1 = 0, \\
\os_{21}^0 = 0, \\
h \si_2 \os_{12}^1 + \tfrac14 h \os_{22}^1 + \os_{22}^0 = 0.
\end{gather*}
Substituting $\os_{12}^0$ and $\os_{22}^0$ and solving the resulting
set of equations yields
$$
\os_{11}^0 = \os_{21}^0 = 0, \quad
\os_{12}^1 = \frac{12(\sqrt3-1)}{h^2} \oy_{k2} \quad\text{and}\quad
\os_{22}^1 = -\frac{12(\sqrt3+1)}{h^2} \oy_{k2}.
$$
\item At order~$\ka^{-1}$, we find
\begin{align*}
& \os_{11}^1 + \tfrac14h\phi_1\os_{11}^0 + \tfrac14h\phi_1\os_{12}^0 +
h\si_1\phi_1\os_{21}^0 + h\si_1\phi_1\os_{22}^0 = -\phi_1\oy_{k1} -
\phi_1\oy_{k2},
\\
& {-}\tfrac14h\phi_1\os_{11}^0 + \tfrac14h\os_{12}^2 + \os_{12}^1 -
\tfrac14h\phi_1\os_{12}^0 - h\si_1\phi_1\os_{21}^0 + h\si_1\os_{22}^2
- h\si_1\phi_1\os_{22}^0 
\\
& \hspace{55ex} = \phi_1\oy_{k1} + \phi_1\oy_{k2},
\\
& h\si_2\phi_2\os_{11}^0 + h\si_2\phi_2\os_{12}^0 + \os_{21}^1 +
\tfrac14h\phi_2\os_{21}^0 + \tfrac14h\phi_2\os_{22}^0 =
-\phi_2\oy_{k1} - \phi_2\oy_{k2},
\\
& {-}h\si_1\phi_2\os_{11}^0 + h\si_2\os_{12}^2 - h\si_2\phi_2\os_{12}^0
- \tfrac14h\phi_2\os_{21}^0 + \tfrac14h\os_{12}^2 + \os_{12}^1 -
\tfrac14h\phi_2\os_{22}^0 
\\
& \hspace{55ex} = \phi_1\oy_{k1} + \phi_1\oy_{k2}.
\end{align*}
We substitute all the known quantities in these equation. The values
of~$\os_{11}^1$ and~$\os_{21}^1$ can then be found from the first and
third equation, respectively: 
$$
\os_{11}^1 = -\phi_1\oy_{k1} \quad\text{and}\quad 
\os_{21}^1 = -\phi_2\oy_{k1}.
$$
The second and fourth equation become
\begin{gather*}
\tfrac14h\os_{12}^2 + h\si_1\os_{22}^2 = \phi_1\oy_{k1} -
12(\sqrt3-1)h^{-2}\oy_{k2}, \\
h\si_2\os_{12}^2 + \tfrac14h\os_{22}^2 = \phi_2\oy_{k1} +
12(\sqrt3+1)h^{-2}\oy_{k2}.
\end{gather*}
The solution of this system is
\begin{align*}
\os_{12}^2 &= \frac1h \Bigl( 3\phi_1 + (2\sqrt3-3)\phi_2 \Bigr)
\oy_{k1} - \frac{24}{h^3} (2\sqrt3-3) \oy_{k2}, \\
\os_{22}^2 &= -\frac1h \Bigl( (2\sqrt3+3)\phi_1 - 3\phi_2 \Bigr)
\oy_{k1} + \frac{24}{h^3} (2\sqrt3+3) \oy_{k2}.
\end{align*}
\item At order~$\ka^{-2}$, the first and third equations are
\begin{gather*}
\os_{11}^2 + \tfrac14h\phi_1\os_{11}^1 + \tfrac14h\phi_1\os_{12}^1 +
h\si_1\phi_1\os_{21}^1 + h\si_1\phi_1\os_{22}^1 = 0,
\\
h\si_2\phi_2\os_{11}^1 + h\si_2\phi_2\os_{12}^1 + \os_{21}^2 +
\tfrac14h\phi_2\os_{21}^1 + \tfrac14h\phi_2\os_{22}^1 = 0.
\end{gather*}
The values of $\os_{11}^2$ and $\os_{21}^1$ follow immediately:
\begin{align*}
\os_{11}^2 &= h \Bigl( \tfrac14\phi_1^2 + (\tfrac14-\tfrac16\sqrt3)
\phi_1\phi_2 \Bigr) \oy_{k1} - \frac{2\sqrt3}{h} \phi_1 \oy_{k2}, \\
\os_{21}^2 &= h \Bigl( (\tfrac14+\tfrac16\sqrt3) \phi_1\phi_2 +
\tfrac14\phi_2^2 \Bigr) \oy_{k1} + \frac{2\sqrt3}{h} \phi_2 \oy_{k2}.
\end{align*}
\end{itemize}
Collecting the results, we find that the stage values for the
Gauss--Legendre method are
\begin{align*}
\os_{11} &= -\frac{\phi_1}{\ka} \oy_{k1} 
+ \frac{h \bigl(3\phi_1^2 + (3-2\sqrt3) \phi_1\phi_2 \bigr)}{6\ka^2} 
\oy_{k1} - \frac{2\sqrt3\,\phi_1}{h\ka^2} \oy_{k2} + \cO(\ka^{-3}),
\\
\os_{12} &= -\frac{2\sqrt3}{h} \oy_{k2} + \frac{12(\sqrt3-1)}{h^2\ka}
\oy_{k2} 
\\
&\hspace{10ex} + \frac{3\phi_1 + (2\sqrt3-3)\phi_2}{h\ka^2} \oy_{k1} -
\frac{24(2\sqrt3-3)}{h^3\ka^2} \oy_{k2} + \cO(\ka^{-3}),
\\
\os_{21} &= -\frac{\phi_2}{\ka} \oy_{k1} 
+ \frac{h \bigl((3+2\sqrt3) \phi_1\phi_2 + 3\phi_2^2 \bigr)}{6\ka^2}
\oy_{k1} + \frac{2\sqrt3\,\phi_2}{h\ka^2} \oy_{k2} + \cO(\ka^{-3}),
\\
\os_{22} &= \frac{2\sqrt3}{h} \oy_{k2} - \frac{12(\sqrt3+1)}{h^2\ka}
\oy_{k2} 
\\
&\hspace{10ex} - \frac{(2\sqrt3+3)\phi_1 - 3\phi_2}{h\ka^2} \oy_{k1} +
\frac{24(2\sqrt3+3)}{h^3\ka^2} \oy_{k2} + \cO(\ka^{-3}).
\end{align*}
The result of doing one step is therefore
\begin{align*}
\oy_{k+1,1} &= \oy_{k1} + \tfrac12h(\os_{11}+\os_{21}) 
\\
&= \oy_{k1} - \frac{h(\phi_1+\phi_2)}{2\ka} \oy_{k1} +
\frac{h^2(\phi_1+\phi_2)^2}{4\ka^2} \oy_{k1} +
\frac{\sqrt3(\phi_2-\phi_1)}{\ka^2} \oy_{k2} + \cO(\ka^{-3})
\\
\intertext{and}
\oy_{k+1,2} &= \oy_{k2} + \tfrac12h(\os_{12}+\os_{22}) 
\\
&= \oy_{k2} - \frac{12}{h\ka} \oy_{k2} +
\frac{\sqrt3(\phi_2-\phi_1)}{\ka^2} \oy_{k1} + \frac{72}{h^2\ka^2}
\oy_{k2} + \cO(\ka^{-3}). 
\end{align*}
We can write this as $\oy_{k+1} = \Psi^-_k \oy_k$ with
$$
\Psi^-_k = \begin{bmatrix} 
  1 - \dfrac{h\al_k}{\ka} + \dfrac{h^2\al_k^2}{2\ka^2} 
  & \dfrac{12\be_k}{h\ka^2} \\[3\jot]
  \dfrac{12\be_k}{h\ka^2} & 1 - \dfrac{12}{h\ka} + \dfrac{72}{h^2\ka^2}
\end{bmatrix} + \cO(\ka^{-3}),
$$
where $\al_k$ and $\be_k$ are defined by~\eqref{albe}. It follows that
$$
\Psi^-_k\oy(\xi_k) = \begin{bmatrix}
  1 - \dfrac{\Phi_-(\xi_k)+h\al_k}{\ka}
  + \dfrac{\bigl(\Phi_-(\xi_k)+h\al_k\bigr)^2}{2\ka^2} \\[3\jot]
  \dfrac{h\phi_-(\xi_k)+12\be_k}{h\ka^2}
\end{bmatrix} + \cO(\ka^{-3}).
$$
Finally, the local error is given by
$$
\oL_k = \Psi^-_k\oy(\xi_k) - \oy(\xi_{k+1}) = \begin{bmatrix} 
  \ka^{-1} \oL_k^{a,-} + \ka^{-2} \oL_k^{b,-} + \cO(\ka^{-3}h^5) \\[\jot]
  \ka^{-2} \oL_k^{c,-} + \cO(\ka^{-3}h^2) 
\end{bmatrix}
$$
where
\begin{align*}
\oL_k^{a,-} &= \Phi_-(\xi_{k+1}) - \Phi_-(\xi_k) - h\al_k 
\\
&= \int_{\xi_k}^{\xi_k+h} \phi(x) \,\d{x} -
\tfrac12h\bigl(\phi(\xi_k^1)+\phi(\xi_k^2)\bigr)
\intertext{and}
\oL_k^{b,-} &= \tfrac12\bigl(\Phi_-(\xi_k)+h\al_k\bigr)^2 
- \tfrac12\bigl(\Phi_-(\xi_{k+1})\bigr)^2 
\\
&= \tfrac12 \Bigl( h\al_k + \Phi_-(\xi_k) - \Phi_-(\xi_{k+1}) \Bigr) 
\Bigl( h\al_k + \Phi_-(\xi_k) + \Phi_-(\xi_{k+1}) \Bigr) 
\\
&= - \oL_k^{a,-} \Bigl( \tfrac12\oL_k^{a,-} + \Phi_-(\xi_k) +
\tfrac12h\bigl(\phi_-(\xi_k^1)+\phi_-(\xi_k^2)\bigr) \Bigr)
\\
\intertext{and}
\oL_k^{c,-} &= \phi_-(\xi_k) + \frac{12\be_k}{h} - \phi_-(\xi_{k+1}) \\
&= \phi(\xi_k) - \phi(\xi_{k+1}) +
\sqrt3\bigl(\phi(\xi_k^1)-\phi(\xi_k^2)\bigr). 
\end{align*}

\subsection{The global error}

The global error is given by the recursion
$$
\oE^-_{k+1} = \Psi^-_k \oE^-_k + \oL_k, \qquad \oE^-_0 = 0.
$$
The solution of this recursion relation is
$$
\oE^-_k = \begin{bmatrix}
\displaystyle \frac1{\ka_-} \sum_{j=0}^{k-1} \oL_j^{a,-} + \frac1{\ka_-^2}
\sum_{j=0}^{k-1} \biggl( \oL_j^{b,-} - h \al_j \sum_{i=0}^{j-1} \oL_i^{a,-}
\biggr) + \cO(\ka_-^{-3}h^4) \\
\displaystyle \frac1{\ka_-^2} \sum_{j=0}^{k-1} \oL_j^{c,-} + \cO(\ka_-^{-3}h^2)
\end{bmatrix}.
$$
Indeed, assuming that the result holds for some value of~$k$, we have
\begin{align*}
\oE^-_{k+1} &= \Psi^-_k \oE^-_k + \oL_k
\\
&= \begin{bmatrix} 
1 - \ka_-^{-1}h\al_k + \cO(\ka_-^{-2}) & \cO(\ka_-^{-2}) \\[\jot]
\cO(\ka_-^{-2}) & 1 + \cO(\ka_-^{-1})
\end{bmatrix}
\\
&\qquad {} \times \begin{bmatrix}
\ka_-^{-1} \sum_{j=0}^{k-1} \oL_j^{a,-} + \ka_-^{-2}
\sum_{j=0}^{k-1} \bigl( \oL_j^{b,-} - h \al_j \sum_{i=0}^{j-1} \oL_i^{a,-}
\bigr) + \cO(\ka_-^{-3}h^4) \\[\jot]
\ka_-^{-2} \sum_{j=0}^{k-1} \oL_j^{c,-} + \cO(\ka_-^{-3}h^2)
\end{bmatrix}
\\
&\qquad {} +  \begin{bmatrix}
\ka_-^{-1} \oL_k^{a,-} + \ka_-^{-2} \oL_k^{b,-} + \cO(\ka_-^{-3}h^5) \\[\jot]
\ka_-^{-2} \oL_k^{c,-} + \cO(\ka_-^{-3}h^2) 
\end{bmatrix}
\\
&= \begin{bmatrix}
\ka_-^{-1} \sum_{j=0}^k \oL_j^{a,-} + \ka_-^{-2}
\sum_{j=0}^k \bigl( \oL_j^{b,-} - h \al_j \sum_{i=0}^{j-1} \oL_i^{a,-}
\bigr) + \cO(\ka_-^{-3}h^4) \\[\jot]
\ka_-^{-2} \sum_{j=0}^k \oL_j^{c,-} + \cO(\ka_-^{-3}h^2)
\end{bmatrix}
\end{align*}
and the formula for the global error follows by induction.

\subsection{The error on $[0,\infty)$}

The solution on the interval $[0,\infty)$ can be computed by running
the Gauss--Legendre method backwards:
\begin{align*}
\os_1 &= \oA_+(\xi_k^1) \, 
\bigl( \oy_k - \tfrac14h \os_1 - (\tfrac14-\tfrac{\sqrt3}6)h \os_2 \bigr), \\
\os_2 &= \oA_+(\xi_k^2) \, 
\bigl( \oy_k - (\tfrac14+\tfrac{\sqrt3}6)h \os_1 - \tfrac14h \os_2 \bigr), \\
\oy_{k+1} &= \oy_k - \tfrac12h (\os_1 + \os_2),
\end{align*}
with $\oA_+$ as given in~\eqref{oAp} and $\xi_k$, $\xi_k^1$ and
$\xi_k^2$ as given in~\eqref{xikp}. The global error can be computed
as before, but here we will take a short-cut. If we comparing the
matrix~$\oA_-$ with~$\oA_+$ and the exact solution on the
interval~$(-\infty,0]$ with the exact solution on~$[0,\infty)$, we
find that they can be related by swapping the components~1 and~2,
replacing $\xi$ by~$-\xi$, and replacing the $-$~subscript with a
$+$~subscript. Hence, the global error of the Gauss--Legendre method
run backwards is
$$
\oE^+_k = \begin{bmatrix}
\displaystyle \frac1{\ka_+^2} \sum_{j=0}^{k-1} \oL_j^{c,+} + \cO(\ka_+^{-3}h^2)
\\
\displaystyle \frac1{\ka_+} \sum_{j=0}^{k-1} \oL_j^{a,+} + \frac1{\ka_+^2}
\sum_{j=0}^{k-1} \biggl( \oL_j^{b,+} - h \al_j \sum_{i=0}^{j-1} \oL_i^{a,+}
\biggr) + \cO(\ka_+^{-3}h^4) 
\end{bmatrix} 
$$
where
\begin{align*}
\oL_k^{a,+} &= \int_{\xi_k-h}^{\xi_k} \phi(x) \,\d{x} -
\tfrac12h\bigl(\phi(\xi_k^1)+\phi(\xi_k^2)\bigr)
\\
\oL_k^{b,+} &= - \oL_k^{a,+} \Bigl( \tfrac12\oL_k^{a,+} + \Phi_-(\xi_k) +
\tfrac12h\bigl(\phi_-(\xi_k^1)+\phi_-(\xi_k^2)\bigr) \Bigr) 
\\
\oL_k^{c,+} &= \phi(\xi_k) - \phi(\xi_{k+1}) +
\sqrt3\bigl(\phi(\xi_k^1)-\phi(\xi_k^2)\bigr). 
\end{align*}

\subsection{The error in the Evans function}

The error in the Evans function is given by~\eqref{ed}:
\begin{equation}
\tag{\ref{ed}}
\begin{aligned}
E_D &= \tfrac12(\ka_- - \ka_+) \Bigl( 
\ov_-(0) \, [\oE^+_k]_2 - \ou_-(0) \, [\oE^+_k]_1 
+ [\oE^-_k]_2 \, \ov_+(0) \\ 
&\hspace{20ex} - [\oE^-_k]_1 \, \ou_+(0) 
+ [\oE^-_k]_2 \, [\oE^+_k]_2 - [\oE^-_k]_1 \, [\oE^+_k]_1
\Bigr) \\
&\quad + \tfrac12(\ka_- + \ka_+) \Bigl( 
\ov_-(0) \, [\oE^+_k]_1 - \ou_-(0) \, [\oE^+_k]_2 
+ [\oE^-_k]_2 \, \ou_+(0) \\ 
&\hspace{20ex} - [\oE^-_k]_1 \, \ov_+(0) 
+ [\oE^-_k]_2 \, [\oE^+_k]_1 - [\oE^-_k]_1 \, [\oE^+_k]_2
\Bigr).
\end{aligned}
\end{equation}
Estimating all the terms, we find that
$$
E_D = - \tfrac12(\ka_- + \ka_+) \Bigl( \ou_-(0) \, [\oE^+_k]_2 
+ [\oE^-_k]_1 \, \ov_+(0) \Bigr) + \cO(\la^{-1/2}h^8, \la^{-3/2}h^2).
$$
Let $\cX$ denote the expression between the big parentheses. We need
to evaluate this expression:
\begin{align*}
\cX &= \ou_-(0) \, [\oE^+_k]_2 + [\oE^-_k]_1 \, \ov_+(0) \\
&= \biggl( 1 - \frac{\Phi_-(0)}{\ka_-} \biggr)
\Biggl( \frac1{\ka_+} \sum_{j=0}^{N-1} \oL_j^{a,+} + \frac1{\ka_+^2}
\sum_{j=0}^{N-1} \biggl( \oL_j^{b,+} - h \al^+_j \sum_{i=0}^{j-1} \oL_i^{a,+}
\biggr) \Biggr) \\
& \qquad + \Biggl( \frac1{\ka_-} \sum_{j=0}^{N-1} \oL_j^{a,-} +
\frac1{\ka_-^2} \sum_{j=0}^{N-1} \biggl( \oL_j^{b,-} - h \al^-_j
\sum_{i=0}^{j-1} \oL_i^{a,-} \biggr) \Biggr) \biggl( 1 -
\frac{\Phi_+(0)}{\ka_+} \biggr) \\
& \qquad + \cO(\la^{-3/2}h^4) \\
&= \frac{1}{\sqrt\la} \sum_{j=0}^{N-1} \Bigl( \oL_j^{a,-} +
\oL_j^{a,+} \Bigr) \\
& \qquad + \frac1\la \sum_{j=0}^{N-1} \Biggl(
\oL_j^{b,-} + \oL_j^{b,+} - h \al^-_j \sum_{i=0}^{j-1} \oL_i^{a,-} 
- h \al^+_j \sum_{i=0}^{j-1} \oL_i^{a,+} \\
& \hspace{30ex} - \Phi_-(0) \oL_j^{a,+} - \Phi_+(0) \oL_j^{a,-}
\Biggr) + \cO(\la^{-3/2}h^4).
\end{align*}
The sum $\sum_j \bigl( \oL_j^{a,-} + \oL_j^{a,+} \bigr)$ is the same as
expression~\eqref{olka}, which appeared in the fourth-order Magnus
method. As we discussed there, this expression is negligible.
Substituting the values of~$\oL_j^{b,-}$ and~$\oL_j^{b,+}$, we find
that
\begin{align*}
\cX &= -\frac1\la \sum_{j=0}^{N-1} \Biggl( \oL_j^{a,-} \Bigl(
\tfrac12\oL_j^{a,-} + \Phi_-(-L+jh) + \al^-_j + \Phi_+(0) \Bigr) +
h \al^-_j \sum_{i=0}^{j-1} \oL_i^{a,-} \\ 
& \hspace{13ex} + \oL_j^{a,+} \Bigl(
\tfrac12\oL_j^{a,+} + \Phi_+(L-jh) + \al^+_j + \Phi_-(0) \Bigr) +
h \al^+_j \sum_{i=0}^{j-1} \oL_i^{a,+} \Biggr)
\end{align*}
Exchanging the double sums yields
\begin{align*}
\cX &= -\frac1\la \Biggl( \sum_{j=0}^{N-1} \oL_j^{a,-} \biggl(
\tfrac12\oL_j^{a,-} + \Phi_-(-L+jh) + \Phi_+(0) + h \sum_{i=j}^{N-1}
\al^-_i \biggr) \\ 
& \hspace{13ex} + \oL_j^{a,+} \biggl( \tfrac12\oL_j^{a,+} +
\Phi_+(L-jh) + \Phi_-(0) + h \sum_{i=j}^{N-1} \al^+_i \biggr) \Biggr).
\end{align*}
Now, $\oL_j^{a,-}$ was defined as
$$
\oL_j^{a,-} = \Phi_-(-L+jh+h) - \Phi_-(-L+jh) - h\al^-_j,
$$
and thus we have
$$
h \sum_{i=j}^{N-1} \al^-_i = \Phi_-(0) - \Phi_-(-L+jh) +
\sum_{i=j}^{N-1} \oL_j^{a,-}.
$$
Using this expression, and its equivalent for $\sum_i \al^-_i$, we
find that
\begin{align*}
\cX &= -\frac1\la \Biggl( \sum_{j=0}^{N-1} \oL_j^{a,-} \biggl(
\Phi_-(0) + \Phi_+(0) + \tfrac12\oL_j^{a,-} + \sum_{i=j}^{N-1}
\oL_i^{a,-} \biggr) \\ 
& \hspace{13ex} + \oL_j^{a,+} \biggl( \Phi_-(0) + \Phi_+(0) +
\tfrac12\oL_j^{a,+} + \sum_{i=j}^{N-1} \oL_i^{a,+} \biggr) \Biggr).
\end{align*}
We know that $\oL_i^{a,\pm} = \cO(h^5)$, which yields
$$
\cX = -\frac1\la \Bigl( \Phi_-(0) + \Phi_+(0) \Bigr) \sum_{j=0}^{N-1}
\Bigl( \oL_j^{a,-} + \oL_j^{a,+} \Bigr) + \cO(\la^{-1}h^8).
$$
The last step is to recall that the sum $\sum_j \bigl( \oL_j^{a,-} +
\oL_j^{a,+} \bigr)$ is negligible, and thus, $\cX = \cO(\la^{-1}h^8)$.
We finally conclude that the error in evaluating the Evans function
with the Gauss--Legendre method is
$$
E_D = \cO(\la^{-1/2}h^8, \la^{-1}h^4, \la^{-3/2}h^2).
$$

\bibliography{../jitse}
\bibliographystyle{habbrv}

\end{document}